\documentclass[secthm,seceqn,amsthm,ussrhead,12pt]{amsart}
\usepackage{amsmath,latexsym}
\usepackage[english]{babel}
\usepackage[psamsfonts]{amssymb}
\usepackage{times}
\usepackage{cite}

\usepackage[mathcal]{euscript}

\newtheorem{Th}{Theorem}
\newtheorem{Lem}[Th]{Lemma}
\newtheorem{proposition}[Th]{Proposition}

\newtheorem{corollary}[Th]{Corollary}

\begin{document}
\sloppy
\thispagestyle{empty}

\textbf{Jordan algebras admitting derivations with invertible values}
\footnote{The work was supported by RFBR 15-31-21169;  FAEPEX 16712-16; FAPESP  14/24519-8, 16/00541-0 and 16/16445-0; and
the President’s Programme ‘Support of Young Russian Scientists’ (grant MK-1378.2017.1).}

\medskip

\medskip
\textbf{Ivan Kaygorodov$^{a}$, Artem Lopatin$^{b,c}$, Yury Popov$^{c,d}$}

\medskip
{\tiny

$^a$ Universidade Federal do ABC, CMCC, Santo Andre, SP, Brazil

$^b$ Omsk Branch of Sobolev Institute of Mathematics, Omsk, Russia

$^c$ Universidade Estadual de Campinas, IMECC, Campinas, SP, Brazil

$^d$ Novosibirsk State University, Novosibirsk, Russia

    E-mail addresses:\smallskip

    Ivan Kaygorodov (kaygorodov.ivan@gmail.com)
    
    Yury Popov (yuri.ppv@gmail.com)

    Artem Lopatin (artem\_lopatin@yahoo.com)
    
}    

\section*{Abstract}

The notion of a derivation with invertible values as a derivation of ring with unity that only takes multiplicatively invertible or zero values appeared in a paper of Bergen, Herstein and Lanski,
in which they determined the structure of associative rings that admit derivations with invertible values.
Later, the results of this paper were generalized in many cases,
for example, for generalized derivations, associative superalgebras, alternative algebras and many others.
The present work is dedicated to the description of all Jordan algebras admitting derivations with invertible values.

\section{Introduction}

Let $A$ be an algebra with unit element $1$ over field $\mathbb{F}.$ 
We denote by $U$ the set of invertible elements of $A$.
Further in this article we only consider \textit{derivations with invertible values},
by which we understand such non-zero derivations $d$ that for every $x \in A$ holds $d(x) \in U$ or $d(x)=0$.

In 1983, Bergen, Herstein and Lanski initiated the study whose purpose is to relate the structure of a ring to the special behavior
of one of its derivations. Namely, in their article \cite{Berg} they described associative rings admitting derivations with invertible values.
They proved that such ring must be either a division ring, or the ring of $2 \times 2$ matrices over a division ring,
or a factor of a polynomial ring over division ring of characteristic 2. They also characterized the division rings
such that the $2 \times 2$ matrix ring over them has an inner derivation with invertible values.
Further, associative rings with derivations with invertible values (and also their generalizations) were discussed in variety of works
(see, for instance, \cite{Car,Chang,CAAF,Hongan,Giam,Komatsu}).
So, in \cite{Giam} semiprime associative rings with involution, allowing a derivation with invertible values
on the set of symmetric elements, were given an examination.
In the work  \cite{Car} Bergen and Carini studied associative rings admitting a derivation with invertible values on some non-central Lie ideal.
Also, in the papers \cite{Chang} and \cite{Hongan} the structure of associative rings that admit $\alpha$-derivations with invertible values and their natural generalizations --- $(\sigma,\tau)$-derivations with invertible values was described.
In the paper \cite{Komatsu} Komatsu and Nakajima described associative rings that allow generalized derivations with invertible values.
The case of associative superalgebras with derivations with invertible values was studied in the paper of Demir, Albas, Argac, and Fosner \cite{CAAF}.
The description of non-associative algebras addmiting derivations with invertible values began in paper of Kaygorodov and Popov \cite{KP}, where it was proved that every alternative (non-associative) algebra addmiting derivation with invertible values
is a Cayley--Dickson over their center or a factor-algebra of polynomial algebra $C[x]/(x^2)$ over a Cayley--Dickson division algebra.

Nowadays, a great interest is shown to the studying of nonassociative algebras and superalgebras with derivations.
For example, in paper \cite{Popov2} the structure of differentiably simple Jordan algebras is determined,
and papers \cite{kay,kay2,shestakov} give the description of generalizations of derivations of Jordan algebras.
Analogues of Moens' theorem, describing nilpotent finite-dimensional algebras as those having invertible Leibniz-derivations were proved for Jordan  algebras \cite{kp15}.
Nevertheless, the problem of specification of Jordan algebras admitting derivations with invertible values remains unconsidered.
The present work is dedicated to the description of Jordan algebras admitting derivations with invertible values.

\section{Preliminaries}

In this article we consider unital algebras over field $\mathbb{F}$ of characteristic $\neq 2.$ Further on by $\bar{\mathbb{F}}$ we denote the algebraic closure of $\mathbb{F}.$
A commutative algebra $J$ is called \textit{Jordan} if it satisfies the \textit{Jordan identity}:
$ (x^2,y,x) = 0.$ All the following information about Jordan algebras in this article can be found in the book \cite{Jac}.

%A derivation $d$ of $A$ is called \textit{inner} if it lies in the smallest subspace of space of all linear operators on $A$,
%containing all right and left multiplications by elements of $A$ and closed under commutation. Otherwise $d$ is called \textit{outer}.

Let $J$ be a unital Jordan algebra. An element $x \in J$ is called \emph{invertible} if there exists $y \in J$ such that $xy = 1, x^2y = x.$

Let $e \in J$ be an idempotent. Then $J$ has the following decomposition:
$$J = J_0 \oplus J_1 \oplus J_2,$$
where $J_i = \{x \in J: ex = \frac{i}{2}x\},$ called the \textit{Peirce decomposition of $J$ with respect to $e$.} The spaces $J_0, J_1$ and $J_2$ satisfy the following multiplicative relations:
\begin{equation}
\label{pd}
J_0J_0 \subseteq J_0, J_2J_2 \subseteq J_2, J_0J_2 = 0, J_0J_1 \subseteq J_1, J_2J_1 \subseteq J_1, J_1J_1 \subseteq J_0 + J_2.
\end{equation}
Analogously, if $J \in 1 = e_1 + \ldots + e_n$ is the sum of $n$ orthogonal idempotents, then $J$ has the following Peirce decomposition with respect to $e_1, \ldots, e_n:$
$$J = \bigoplus_{i,j = 1}^n J_{ij},$$
where
$$J_{ii} = \{x \in J: e_ix = 0, e_jx = 0, j \neq i\},$$
$$J_{ij} = \{x \in J: e_ix = e_jx = \frac{1}{2}x, e_kx = 0, k \neq i, j\} = J_{ji}, i \neq j.$$
%As above, the following inclusions hold:
%\begin{equation}
%\label{pdmult_1}
%J_{ii}^2 \subseteq J_{ii}, J_{ii}J_{jj} = 0, %J_{ii}J_{ij} \subseteq J_{ij}, i \ neq j;
%\end{equation}
%\begin{equation}
%\label{pdmult_2}
%\end{equation}

Let $A$ be a unital algebra over $\mathbb{F}$. The \textit{degree} of $A$ is the maximal number of orthogonal idempotents that sum to 1 in a scalar extension $\bar{A} = A \otimes_F \bar{\mathbb{F}}$.

During our discussion we will encounter some certain types of Jordan algebras, so, in order to make our work self-contained, we provide their definitions:

\begin{enumerate}
\item[I.]  \textbf{Algebras of type $A^{(+)}$.} Let $(A,\cdot)$ be an algebra. Then its underlying vector space equipped with new multiplication $a \circ b = \frac{1}{2}(a\cdot b + b \cdot a)$ is again an algebra which is denoted by $A^{(+)}.$ If $A$ is associative, then $A^{(+)}$ is a Jordan algebra.
If a Jordan algebra $J$ can be embedded into $A^{(+)}$ for an associative algebra $A$, it is called \textit{special}. Nonspecial Jordan algeras are called \textit{exceptional}. It is known that an element of a special Jordan algebra $J \subseteq A^{(+)}$ is invertible if and only if it is invertible in $A$ with the same inverse. It is known that if $A$ is unital associative, the degree of $A^{(+)}$ is the same that of $A.$
%The subalgebra generated by set $J$ in $A$ is called \textit{associative enveloping algebra} for J.

\item[II.] \textbf{Algebras of type $H(A,*)$.} For an algebra $A$ with involution $*$ the set $H(A,*)$ of $*$-hermitian elements
(i.e., such elements $x$ that $x^* = x$) is closed under Jordan product $\circ$, so it is a subalgebra of $A^{(+)}$. Again, if $A$ is associative, $H(A,*)$ is a Jordan algebra.

During the paper we will frequently encounter matrix algebras with involutions, so it is reasonable to give some examples of them:

1) The first example is the well-known transpose involution $T$. It is obvious that the degree of the algebra $H(M_n,T)$ is exactly $n.$

2) For even values of $n$, the algebra $M_n = M_{2m}$ admits the \textit{symplectic involution} $Symp$, which acts as follows:
$$\begin{pmatrix}A & B \\ C & D \end{pmatrix}^{Symp} = \begin{pmatrix}D^T & -B^T \\ -C^T & A^T \end{pmatrix},$$
where $A, B, C, D \in M_m.$
The degree of the algebra $H(M_{2n},Symp)$ is exactly $n.$

It is known that every finite-dimensional simple associative algebra with involution over an algebraically closed field is isomorphic to either $(M_n,T)$ or to $(M_{2n},Symp)$. Now let $A$ be a simple finite-dimensional algebra over $\mathbb{F}$ with involution $*$. Then $A \otimes_\mathbb{F} \Omega$ is isomorphic either to $(M_n,T)$ or to $(M_{2n},Symp)$. Accordingly, the involution $*$ is called to be of \textit{transposition} or of \textit{symplectic} type.
\item[III.] \textbf{Algebras of symmetric bilinear vector form $J(V,f)$.}
Let $V$ be a $\mathbb{F}$-vector space with a symmetric bilinear form $f:V\otimes V \rightarrow \mathbb{F}$.
We can endow $\mathbb{F} \oplus V$ with a structure of Jordan algebra by defining multiplication as follows:
$$(\alpha + v)(\beta + u) = \alpha\beta + f(v,u) + \alpha u + \beta v, \text{ where } \alpha, \beta \in \mathbb{F}, u, v \in V.$$
We denote this algebra by $J(V,f).$ The degree of the algebra $J(V,f)$ is 2.

\item[IV.] \textbf{Algebras of Albert type.}
Let $C$ be a Cayley-Dickson algebra over $\mathbb{F}$, and $\gamma_1, \gamma_2, \gamma_3$ be nonzero elements of $\mathbb{F}$. We denote the main involution in $C$ by \ $\bar{}$ \ (For more information on Cayley-Dickson algebras, we refer the reader to \cite{Jac, kolca}). Consider $C_3$, the algebra of $3\times 3$ matrices with entries in $C.$ By $\overline{X}^T$ we denote the matrix obtained from $X \in C_3$ by applying the transpose and involution  \ $\bar{}$ \ to every coefficient of $X,$ and by $\gamma$ we denote the matrix $diag\{\gamma_1, \gamma_2, \gamma_3\}.$ Then the mapping $*_\gamma:C\rightarrow C$ defined by $X^{*_\gamma} = \gamma^{-1} \overline{X}^T \gamma $ is an involution of $C_3$. It is known that $H(C_3,*_\gamma)$ is an exceptional simple Jordan algebra. 
If $\gamma_1 = \gamma_2 = \gamma_3 = 1$, and $C$ is a split Cayley-Dickson algebra over $\mathbb{F}$, then this algebra is called the Albert algebra over $\mathbb{F}.$

 A Jordan algebra $J$ over $\mathbb{F}$ is called an \textit{algebra of Albert type}, if a scalar extension $J \otimes_{\mathbb{F}} \bar{\mathbb{F}}$ is an Albert algebra over $\bar{\mathbb{F}}$.
 It is known that all algebras of Albert type are divided in two classes: the division algebras and the algebras $H(C_3, *_\gamma)$.

 An element $A \in J = H(C_3,*_\gamma)$ has the form 
$$ A = \begin{pmatrix}
\alpha_1 & c & \gamma_1^{-1}\gamma_3\bar{b} \\
\gamma_2^{-1}\gamma_1\bar{c} & \alpha_2 & a \\
b & \gamma_3^{-1}\gamma_1\bar{a} & \alpha_3 \end{pmatrix}, $$
where $\alpha_1, \alpha_2, \alpha_3 \in F, a, b, c \in C.$
One can see that this algebra is of degree 3, that is, it has 3 orthogonal idempotents $e_{ii}, i = 1, 2, 3$. Relative to these idempotents it has the Peirce decomposition
$$J = \sum_{i=1}^3 J_{ii} + \sum_{i<j} J_{ij},$$
where $J_{ii} = Fe_{ii},$ and $J_{ij} = \{ae_{ij} + \gamma_i\gamma_j^{-1}ae_{ji}, a \in C\}, i, j = 1, 2, 3.$

For an element $A \in J$ we can define
$$n(A) = \alpha_1\alpha_2\alpha_3 - \alpha_1\gamma_3^{-1}\gamma_2n(a) - \alpha_2\gamma_1^{-1}\gamma_3n(b) - \alpha_3\gamma_2^{-1}\gamma_1n(c) + t((ca)b),$$
where for $n(x), x = a, b, c$ we mean the norm of an element of a Cayley-Dickson algebra.
It is well known that $A$ is invertible in $J$ if and only if $n(A) \neq 0.$
\end{enumerate}

We will also need the following statement which describes simple Jordan algebras:\\

{\bf Zelmanov's Theorem \cite{Zel83}.}
Let $J$ be a simple Jordan algebra over field of characteristic $\neq 2$. Then one of the following holds:

\begin{enumerate}
\item $J$ is an algebra $A^{(+)}$, where $A$ is a simple associative algebra;
\item $J$ is an algebra $H(A, *)$, where $A$ is a simple associative algebra;
\item $J$ is an algebra of non-degenerate symmetric bilinear form $J(V,f)$ on a vector space $V$ of dimension $>1$ over field $\mathbb{F}$;
\item $J$ is an algebra of Albert type.
\end{enumerate}

\section{Jordan algebras addmiting derivations with invertible values}

The purpose of this paragraph is to generalize the results of Bergen, Herstein and Lanski to the Jordan case. Further on, by $J$ we denote a Jordan algebra of characteristic $\neq 2$ with unit element $1$ and a derivation with invertible values $d$.

Now we shall study the ideal structure of $J$:
%\begin{Lem}
%\label{Lem1.2}
%Let $I\lhd J$. Suppose that $I \not \subset \ker (d)$. Then:
%\begin{enumerate}
%\item[(a)] $I$ is both minimal and maximal,
%\item[(b)] $I^3=0$.
%\end{enumerate}

%\end{Lem}

%\textbf{Proof.}
%$(a).$
%Let $I_1 \subseteq I \subseteq I_2$ be proper ideals in $J$, and $d(I) \neq 0$. Since $I$ does not contain invertible elements, it is easy to see that $d(I) \cap I=0$ and $I \oplus d(I)$ is also an ideal in $J$.
%Since $d(I) \neq 0$, $d(I)$ contains invertible elements, hence  $I \oplus d(I)=J$. Particularly, for any $j \in J$ there exist unique $a, b \in I$ such that $d(j)=a+d(b)$. Hence $a=d(j-b) \in I \cap d(J)=0$, therefore $J=I \oplus \ker (d)$ which implies $I \cap \ker (d) = 0$.
%For arbitrary $j \in I_2 $ we have $j=m+d(n)$, $m, n \in I$. Consequently, $d(n)=j-m \in I_2 \cap d(I) = 0$; thus $j=m \in I$, so $I$ is maximal.
%Since $I \cap \ker (d) = 0$, $I_1 \not\subset \ker (d)$, therefore, repeating the argument above, we conclude that $I_1 = I$, so $I$ is minimal.

%$(b).$ The cube of an ideal in Jordan algebra is also an ideal.
 %Obviously, $I^3 \subseteq I$, and $I^3$ is also an ideal of $J$, so by $(a)$ we only have two possibilities: $I^3 = I$ or $I^3 = 0$. If $I^3 = I$, then
%$$d(I) = d(I^3) = d(I^2I) \subseteq I^2d(I) + (d(I)I)I \subseteq I \cap d(J) = 0,$$
%which contradicts the initial condition $I \not\subset ker(d)$, so we are forced to conclude that $I^3 = 0$. \hfill$\Box$\vspace{0.3cm}

\medskip

\begin{Lem}
\label{id_ker}
%Suppose $char(F) \neq 2,3$. Then
$d(I) = 0$ for any proper ideal $I$ of $J.$
\end{Lem}

\textbf{Proof.}
%Let $I \lhd J$ and suppose that $d(I) \neq 0.$ Then there exists $i \in I$ such that $d(i) \in U$. Since $I^3 = 0$,
%\begin{eqnarray*}0=d^3(i^3)=
%d^3(i^2)i + 6 d(i)^3 + 6(id^2(i))d(i) + 6(id(i))d^2(i) + i^2 d^3(i), \end{eqnarray*}
%and we have $6d(i)^3 \in I$. Since the power of an invertible element in Jordan algebra is also invertible, 
%$6d(i)^3 = 0$, therefore $6 = 0.$ We have obtained a contradiction which proves the lemma. \hfill$\Box$\vspace{0.3cm}
Let $x \in I.$ Since $I$ is proper, it does not contain invertible elements of $J.$ Then $d(x^2) = 2xd(x) \in I \cap d(J) = 0.$ Hence $0 = d(xd(x)) = d(x)^2 + xd^2(x),$ therefore $d(x)^2 \in I.$ Since the power of an invertible element of a Jordan algebra is also invertible, we must have $d(x) = 0.$ The lemma is now proved.

\begin{Lem}
\label{id_square}
$I^2 = 0$ for any proper ideal $I$ of $J.$
\end{Lem}

\textbf{Proof.}
For any $x \in I, y \notin \ker (d)$ we have $0 = d(xy) = xd(y)$. Hence, we have $xd(y) = x^2d(y) = 0$.
It follows that $(d(y), J, x^2) = 0$  (see \cite{Zel78}).
Particularly, we have
 $$0 = d(y)(d(y)^{-1}x^2) = (d(y)d(y)^{-1})x^2 = x^2.$$
Linearizing this, we obtain $2ab = 0$ for $a, b \in I$, so $I^2 = 0$.\hfill$\Box$\vspace{0.3cm}

\medskip

%We now consider the case when $char(F) \neq 3$.
Denote by $M$ the sum of all ideals contained in $\ker (d)$.
Lemma \ref{id_ker} implies that $M$ is the largest ideal of $J$.

Since $d(M) = 0$, the map $\bar{d}$ on $\bar{J}$ given by $\bar{d}(j + M) = d(j) + M$ is correctly defined, and it is easy to see that
$\bar{d}$ is a derivation with invertible values of $\bar{J}$. Since $M$ is the largest ideal of $J,$ $\bar{J} = J/M$ is a simple Jordan algebra that admits a derivation with invertible values.

We have proved our first result, which is simple but rather important:
\begin{Th}
\label{first_th}
Let $J$ be a Jordan algebra admitting a derivation with invertible values $d$. Then $J$ is an extension of a simple Jordan algebra with derivation with invertible values by an ideal $M$ such that $M^2 = 0, d(M) = 0$ and $M$ is the largest ideal of $J.$
\end{Th}
\medskip

\section{Simple Jordan algebras admitting derivations with invertible values}
Further in this section we consider only simple Jordan algebras admitting derivations with invertible values.

As Zelmanov's theorem suggests, we have to study four cases:

\subsection{The case $J = A^{(+)}$, $A$ is a simple associative algebra} This case is the easiest one, because
by \cite{Herstein} every derivation of $A^{(+)}$ is a derivation of $A$. Since any element which is invertible in $A^{(+)}$ is also invertible in $A$, then $d$ is a derivation with invertible values of associative algebra $A$, so, by \cite{Berg} $A$ is either a division algebra $D$, or a $D_2$, the $2 \times 2$ matrix algebra over a division algebra $D$.

\subsection{The case $J = H(A, *)$, $A$ is a simple associative algebra} This case is dealt with in
\begin{Lem}
\label{herm_der}
Let $A$ be a simple associative algebra such that $H(A, *)$ has a derivation $d$ with invertible values. Then $A$ is isomorphic either to
\begin{enumerate}
\item  a division algebra $D$, or
\item $D_2$, $2 \times 2$ matrix algebra over $D$, or
\item $Z_4$, $4 \times 4$ matrix algebra over its center $Z(A)$, $*$ is the involution of symplectic kind. Moreover, in this case $H(A,*)$ is isomorphic to an algebra of a symmetric nondegenerate bilinear form.
\end{enumerate}
\end{Lem}

\textbf{Proof.}
A derivation $d$ of $H(A, *)$ can be considered as a Jordan derivation from $H(A, *)$ to $A$. By \cite{Lagutina} either any Jordan derivation $d:H(A, *) \rightarrow A$ can be extended to derivation $d:alg(H(A, *)) \rightarrow A$, (by $alg(H(A, *))$ we mean the subalgebra of $A$ generated by $H(A, *)$) or $A$ is a central order in $(Z(Z^{\times})^{-1})_4$, where $Z$ is the center of $A$ and $Z^{\times} = Z \setminus \{0\}$ (that is, $A(Z^{\times})^{-1}$ is isomorphic to the $4 \times 4$ matrix algebra over $Z(Z^{\times})^{-1}$) and $*$ is an involution of symplectic kind. But in the second case, since $A$ is simple, $Z$ is a field, hence $A$ is isomorphic to the $4 \times 4$ matrix algebra over $Z.$ Since $*$ is of symplectic kind, the scalar extension of $A$ by the algebraic closure $\bar{Z}$ of $Z$ is isomorphic to $(\bar{Z}_4,Symp).$ It is easy to see that $H(\bar{Z}_4,Symp)$ is an algebra of symmetric bilinear form on the space of $Symp$-symmetric traceless matrices. Therefore (see \cite{Jac}) $H(A,*)$ is also an algebra of symmetric bilinear form. Thus, further in this lemma we suppose that $d$ extends to a derivation $d:alg(H(A, *)) \rightarrow A$. From Herstein's work \cite{Herstein1} it follows that if $dim_Z A >4$, then $alg(H(A, *)) = A$. Wedderburn--Artin theorem implies that if $dim_Z A \leq 4$, then $A$ is either a division algebra over $Z$ or $Z_2$, which correspondingly matches the cases $(1)$ and $(2)$, so from now on we may assume that $d$ can be extended to a derivation of $A$. As a derivation of $A$, $d$ has invertible values on symmetric elements, so by \cite{Giam} $A$ is an algebra of type $(1)$, $(2)$ or $(3)$.\hfill$\Box$\vspace{0.3cm}

\subsection{The case $J = J(V,f)$} In this section we find necessary and sufficient condintions for the algebra $J$ which is an algebra of nondegenerate bilinear form over its center $Z$ to have a derivation with invertible values, and also characterize all derivations with invertible values of this algebra in the finite-dimensional case.

First we need some preparatory work. Let $J$ be an algebra of symmetric nondegenerate bilinear form $f$ on a space $V$ over its center $Z$ (which is an extension of $\mathbb{F}$).
Define the form $f$ on the whole $J$. Particularly, we set $f(Z,V) = 0,$ and $f(z,z') = -zz'$ for $z, z' \in Z.$ It is easy to see $f$ is still nondegenerate. For an element $x = z + v, z \in Z, v \in V$ define
\begin{equation}
\label{trace}
t(x) = 2z = -2f(x,1),
\end{equation}
\begin{equation}
\label{norm}
n(x) = z^2 - f(v,v) = -f(x,x).
\end{equation}
The elements $t(x)$ and $n(x)$ lie in $Z$ and are called respectively the \textit{trace} and the \textit{norm} of $x.$ The function $t$ is $Z$-linear, the function $n$ is a quadratic form over $Z$ and the element $x$ is invertible if and only if $n(x) \neq 0.$ 
Moreover, one can check that the following relation holds in $J:$
\begin{equation}
\label{sqr_eqn}
x^2 - t(x)x + n(x) = 0.
\end{equation}
The relation (\ref{norm}) implies that
\begin{equation}
\label{f_norm}
n(x+y) - n(x) - n(y) = -2f(x,y).
\end{equation}

The linearized version of relation (\ref{sqr_eqn}) is the following:
\begin{equation}
\label{f_norm_linearized}
xy +yx - t(x)y - t(y)x - 2f(x,y) = 0.
\end{equation}
Let $y = d(x)$ in the relation (\ref{f_norm_linearized}):
\begin{equation}
\label{sqr_eqn_lin_xd(x)}
2xd(x)-t(x)d(x) - t(d(x))x - 2f(x,d(x)) = 0.
\end{equation}
Apply $d$ to the relation (\ref{sqr_eqn}):
\begin{equation}
\label{sqr_eqn_d}
2xd(x) -d(t(x))x - t(x)d(x) + d(n(x)) = 0.
\end{equation}
Subtracting (\ref{sqr_eqn_d}) from (\ref{sqr_eqn_lin_xd(x)}), we get 
$$(d(t(x)) - t(d(x)))x -2f(x,d(x)) - d(n(x)) = 0.$$
If $x \in Z,$ then $d(t(x)) = t(d(x)), -2f(x,d(x)) = d(n(x)).$ If $x \in Z,$ then $t(x) = 2x, n(x) = x^2$ and these relations are trivial. Thus, we have 
\begin{equation}
\label{dt}
t(d(x)) = d(t(x)),
\end{equation}
\begin{equation}
\label{dn}
d(n(x)) = -2f(x,d(x)).
\end{equation}
%Using the relation (\ref{f_norm}), linearize the last relation:
%\begin{equation}
%\label{df}
%d(f(x,y)) = f(x,d(y)) + f(d(x),y).
%\end{equation}
The first result of the section concerns the linearity of derivations with invertible values over the center of an algebra. 

\begin{Lem}
\label{form_der_lin}
Let $J$ be an algebra $J(V,f)$ of nondegerate bilinear form $f$ over its center $Z$, and $d$ be a derivation with invertible values of $J.$ Suppose that dim$_Z J = 2n < \infty$ and $J$ contains an $f$-isotropic space $V$ of dimension $n$ (in other words, $f$ is of maximal Witt index). Then $d$ is linear over $Z$.
\end{Lem}
\textbf{Proof.} Note that for $x \in J, z \in Z$ we have $d(zx) - zd(x) = d(z)x.$ Hence we must prove that $d(Z) = 0.$
Note that $d([a,b]) = [d(a),b] + [a,d(b)], ~ d((a,b,c)) = (d(a),b,c) + (a,d(b),c) + (a,b,d(c)),$ where $a, b, c \in J, [a,b] = ab-ba$ is the \textit{commutator} of $a, b;$ $(a,b,c) = (ab)c - a(bc)$ is the \textit{associator} of $a, b, c.$ Using this formulas, it is easy to see that $d(Z) \subseteq Z.$

Suppose that there exists $0 \neq x \in \ker d = B$ which is not invertible and let $z \in Z.$ Then $d(zx) = d(z)x$ is not invertible, hence it is equal to 0, thus $d(z) = 0.$ The arbitrarity of $z$ implies that $d(Z) = 0,$ and we are done. Hence we may assume that $B$ is a division subalgebra of $J$. %In particular, this means that the restriction of $f$ to $B$ is nondegenerate. Indeed, if $f(x,B) = 0$ for $x \in B$, then $0 = f(x,x) = \text{(by (\ref{norm}))} = -n(x),$ hence $x$ is not invertible. 

Let $V$ be a $n$-dimensional $Z$-subspace of $J$ without invertible elements. Since $V$ does not contain invertible elements, we get 
\begin{equation}
\label{vcapb}
V \cap B = 0.
\end{equation}

Note that $V + d(V)$ is a $Z$-subspace of $J$. Indeed, for $z \in Z, v \in V, zd(v) = d(zv) - d(z)v$ lies in $V + d(V)$ since $d(z) \in Z, zv, d(z)v \in V.$ Now let $e_1, \ldots e_n$ be a basis of $V.$ We prove that $e_1, \ldots, e_n, d(e_1), \ldots, d(e_n)$ are linearly independent over $Z.$ Indeed, if for $z_1, \ldots, z_n, z_1', \ldots z_n' \in Z$ we have $z_1e_1 + \ldots + z_ne_n + z_1'd(e_1) + \ldots + z_n'd(e_n) = 0$ then 
$$0 = (z_1 - d(z_1'))e_1 + \ldots + (z_n - d(z_n)')e_n + d(z_1'e_1 + \ldots + z_n'e_n).$$
The first $n$ summands lie in $V,$ and the last summand lies in $d(V)$. Since $V$ does not contain invertible elements, $V \cap d(V) = 0,$ hence $d(z_1'e_1 + \ldots + z_n'e_n) = 0.$ Then (\ref{vcapb}) implies that $z_1' = \ldots = z_n' = 0.$ Therefore $z_1e_1 + \ldots + z_ne_n = 0$ and $z_1 = \ldots = z_n = 0.$ Since $\dim_Z J = 2n$, we conclude $V \oplus d(V) = J$ (as a sum of $\mathbb{F}$-subspaces). Therefore for any $x \in J$ there exist $a, b \in V$ such that $d(x) = a + d(b),$ thus $a = d(x-b) \in V \cap d(J) = 0$ and $d(x) = d(b).$ Relation (\ref{vcapb}) implies that such $b$ is unique, hence $J = V \oplus B$ as a sum of $\mathbb{F}$-subspaces. Let $z \in Z$ be such that $d(z) \neq 0.$ Then $z = v + b,$ where $v \in V, b \in B$. Since $v$ is not invertible, we have  
$$0 = d(n(v)) = \text{(by (\ref{dn}))} = -2f(v,d(v)) = -2f(v,d(z)) = \text{(since $d(z) \in Z$)} = -2d(z)f(1,v) = d(z)t(v).$$
Since $d(z) \neq 0,$ we must have $t(v) = 0.$ Therefore $2z = t(z) = t(b).$ Apply $d$ to this relation: $2d(z) = d(t(b)) = \text{(by (\ref{dt}))} = t(d(b)) = 0,$ hence $d(z)  = 0,$ a contradiction. Therefore $d(Z) = 0.$ \hfill$\Box$\vspace{0.3cm}

From now on and to the end of the subsection we consider only the derivations linear over the base field.

\begin{Lem}
\label{form_der_char}
A mapping $d$ is a derivation with invertible values of $J(V,f)$ if and only if\\
$i)$ $d(V) \subseteq V;$\\
$ii)$ $d(F) = 0;$\\
$iii)$ $d$ is skew-symmetric with respect to $f;$\\
$iv)$ $f(v,v) \neq 0$ for all nonzero $v \in Im(d)$.
\end{Lem}
\textbf{Proof}. Let $d$ be a derivation of $J$. Then obviously $d(F)=0.$ From relations (\ref{dt}) and (\ref{dn}) it follows that $d(V) \subseteq V$ and that $d$ is skew-symmetric with respect to $f$. Conversely, it is easy to check that
for arbitrary endomorphism $d$ of $J$ the conditions i) --- iii) imply that $d$ is a derivation of $J(V,f)$.
Moreover, we already noted that $u \in V$ is invertible in $J(V,f)$ if and only if $f(u,u) \neq 0$.
Hence derivations with invertible values of $J(V,f)$ are the maps $d \in End(V)$ such that
\begin{eqnarray}
\label{f1}
f(v,d(v)) = 0,
\end{eqnarray}
\begin{eqnarray}
\label{f2}
f(d(v),d(v)) = 0 \Rightarrow d(v) = 0
\end{eqnarray} for all $v \in V.$ \hfill$\Box$\vspace{0.3cm}
\begin{Lem}
\label{form_der}
An algebra $J(V,f)$ admits a derivation with invertible values if and only if there exist
$x, y \in V$ such that $f(x,x) \neq 0,\ f(y,y) \neq 0,\ f(x,y) = 0$ and $-\frac{f(y,y)}{f(x,x)}$ is not a square in the base field $\mathbb{F}$.
\end{Lem}

\textbf{Proof.}
%Suppose, on the contrary, that $\mathbb{F}$ allows square root extraction.
Suppose that the lemma condition does not hold in $V$. In this case we prove that $dim(d(J)) < 2$. Let $x, y$ be two linearly independent vectors in $d(J)$.
By our hypothesis, $f(x,x) \neq 0$, thus we may substitute $y \rightarrow y - \frac{f(x,y)}{f(x,x)}x$ and assume that $f(x,y) =0$.
Also $f(y,y)$ remains nonzero.
%For $\alpha \in F$ we have $f(\alpha x + y, \alpha x + y) = \alpha^2f(x,x) + 2\alpha f(x,y) + f(y,y)$. Consider this as the quadratic polynomial in variable $\alpha$. By our assumption this polynomial has roots in $\mathbb{F}$, which contradicts (\ref{f2}),so we conclude that $d(J) = Fu$ for $u \in V.$
For $\alpha \in F$ we have $f(\alpha x + y, \alpha x + y) = \alpha^2 f(x,x) + f(y,y)$. But for $\alpha = \sqrt{-\frac{f(y,y)}{f(x,x)} } \in F$
this expression is equal to $0$, hence $\alpha x + y$ is not invertible in $J(V,f)$, which contradicts (\ref{f2}), so we conclude that
$d(J) = Fu$ for $u \in V.$
Particularly, $d(u) = \delta u,$ where $\delta \in F.$ But from (\ref{f1}) it follows that $0 = f(u,d(u)) = \delta f(u,u)$, and since
$f(u,u) \neq 0$, we have $\delta  = 0$. Now, take  $v \in V$ such that $d(v) = u$. Linearizing (\ref{f1}), we have
$0 = d(f(v,u))= f(d(v),u) + f(v,d(u)) = f(u,u)$, a contradiction.

Conversely, suppose that $x, y \in V$ satisfy the lemma condition.
By $W$ we denote the subspace spanned by $x$ and $y$.
Since $W$ is finite-dimensional and the restriction of $f$ to $W$ is nondegenerate, it is easy to see that $W \cap W^{\bot} = 0$
and $V = W \oplus W^{\bot}$ for $W^{\bot}$ the orthogonal complement of $W$ with respect to $f.$

Now we are able to explicitly construct a derivation with invertible values of $J(V,f)$: define the mapping $d$
by $d(x) = y, d(y) = -\frac{f(y,y)}{f(x,x)}x$, $d(W^{\bot}) = 0$.
It is easy to see that conditions (\ref{f1}) and (\ref{f2}) hold for $d$, so it is a derivation of $J(V,f)$.
Also, the lemma condition implies that for $\beta, \gamma  \neq 0$,
$$f(\beta x + \gamma y, \beta x + \gamma y) =  \beta^2 f(x,x) +  \gamma^2  f(y,y)\neq 0.$$
Hence $d$ is a derivation with invertible values of $J$.
\hfill$\Box$\vspace{0.3cm}

\begin{corollary}
Every simple Jordan algebra of symmetric bilinear vector form over perfect (or algebraically closed) field
does not have a derivation with invertible values.
\end{corollary}

%In the next lemma we describe derivations with invertible values of $J(V,f)$ in the case when $V$ is finite-dimensional and $f$ is nondgenerate.

\begin{Lem}
\label{form_findim_der}
Let $J = J(V,f)$ be an algebra of a symmetric nondegenerate bilinear form on a finite-dimensional vector space $V$. Then a subspace $U \subseteq V$ is an image of a derivation with invertible values if and only if dim$U$ is even and $f(x,x) \neq 0$ for any $x \in U.$ 
\end{Lem}
\textbf{Proof.} In the lemma \ref{form_der_char} we noticed that a derivation of $J$ is a skew-symmetric mapping on $V$ with respect to $f.$ Therefore, its image must be even-dimensional. Moreover, from (\ref{f2}) it follows that $f(x,x) \neq 0$ for all $x \in U.$

On the other hand, let $U$ be an even-dimensional subspace of $V$ such that $f(x,x) \neq 0$ for any $x \in U.$ Let $f_U$ denote the restriction of $f$ to $U$. We construct a nondegenerate linear mapping $d$ on $U$ which is skew-symmetric with respect to $f_U.$ First of all, note that the conditions imposed on $U$ imply that $f_U$ is nondegenerate. Hence, one can pick a basis $v_1,\ldots,v_{2k}$ of $U$ such that $f_U(v_i,v_j) = \delta_{ij}\alpha_i,$ where $\alpha_i \neq 0$ for all $i$. Now, one can define $d$ by $d(v_{2i-1}) = v_{2i}, d(v_{2i}) = -\frac{\alpha_i}{\alpha_j}v_{2i-1}, i = 1, \ldots,k.$ Since $f_U$ is nondegenerate, $V = U \oplus U^{\bot},$ and we can extend $d$ on $V$ by setting $d(U^{\bot}) = 0$. It is easy to see that $d$ is skew-symmetric with respect to $f$ and that Im$(d) = U,$ thus $d$ satisfies (\ref{f1}) and (\ref{f2}) and is a derivation with invertible values. \hfill$\Box$\vspace{0.3cm}

\subsection{The case $J$ is an algebra of Albert type}
 As we have noted, the algebras of Albert type are divided in two classes: the division algebras (which are not interesting to us, because all their derivations have invertible values), and the algebras $H(C_3,*_\gamma)$. We deal with the later ones in the 

\begin{Lem}
\label{alb_der}
Let $J$ be an algebra $H(C_3,*_\gamma)$ over its center $Z$. Then $J$ does not have nonzero derivations with invertible values.
\end{Lem}
\textbf{Proof.} Let $C$ be a Cayley-Dickson algebra over $Z,$ and $d$ be a derivation with invertible values of $J.$ Recall that $J$ has the Peirce decomposition $J = \sum_{i=1}^3 J_{ii} + \sum_{i<j} J_{ij}$. Then $d(e_{11}) = d(e_{11}^2) = 2e_{11}d(e_{11}),$ which means that $d(e_{11}) \in J_{12} + J_{13}, $ that is, it has the form
$$d(e_{11}) = \begin{pmatrix}
0 & c & \gamma_1^{-1}\gamma_3\bar{b} \\
\gamma_2^{-1}\gamma_1\bar{c} & 0 & 0 \\
b & 0 & 0 \end{pmatrix}, $$
where $ \gamma_1, \gamma_2, \gamma_3 \in Z; b, c \in C.$

It is easy to see that the norm $n(d(e_{11}))$ as an element of $H(C_3,*_\gamma)$ is zero, therefore, $d(e_{11})$ is not invertible and must be 0. Analogously, $d(e_{22}) = d(e_{33}) = 0.$ Now, let $a \in J_{ii}.$ Then, since $d(e_{ii}) = 0,$ we have $d(a) = d(e_{ii}a) = e_{ii}d(a).$ Hence $d(a) \in J_{ii}.$ Again, the norm of this element is zero, hence it must be zero. We have proved that $d(J_{ii}) = 0, i = 1, 2, 3.$ Analogously, let $a \in J_{ij}, i \neq j.$ Then $d(a) = 2e_{ii}d(a) = 2e_{jj}d(a),$ therefore, $d(a) \in J_{ij}.$ One can see that $n(d(a))$ is again zero, therefore, $d(a)$ is zero. Hence $d(J_{ij}) = 0$ for all Peirce components of $J$, and $d = 0.$ \hfill$\Box$\vspace{0.3cm}

Now we can state the main result of the section:

\begin{Th}
\label{mainth}
Let $J$ be a Jordan algebra of characteristic $\neq 2$ admitting derivation with invertible values $d$.
Then one of the following holds:
\begin{enumerate}
\item  $J$ is an algebra $A^{(+)}$, where $A = D$ or $D_2$, $D$ is an associative division algebra;
\item  $J$ is an algebra $H(A, *)$, where $A = D$ or $D_2$, $D$ is an associative division algebra;
\item  $J$ is an algebra of symmetric nondegenerate bilinear form $J(V,f)$;
\item  $J$ is a division algebra of Albert type;
\item  $J$ is an extension of cases $(1)$ -- $(4)$ by $M = \mathbf{P}(J)$ -- the prime radical of $J$, $M \subseteq \ker d$, $M$ is the largest ideal of $J$.
\end{enumerate}

\end{Th}

\medskip

\textbf{Proof.}
Follows immediately from lemmas \ref{id_ker} -- \ref{alb_der}.
\medskip

\section{The finite-dimensional case}

Now we return to the general situation. Recall that by theorem \ref{first_th} every Jordan algebra $J$ with a derivation with invertible values $d$ is an extension of a simple Jordan algebra
$\bar{J}$ with a derivation with invertible values $\bar{d}$ by a trivial ideal $M.$ In this section we study this extensions assuming the finite-dimensionality of all spaces involved.

We will need some additional terminology.
Following to the general conception due to Eilenberg \cite{Eilenberg}, define  a notion of a Jordan bimodule.
 Let $A$ be a Jordan algebra, and $N$ be a bimodule over $A$. 
We call the algebra $E = A \oplus N$  with the multiplication extending multiplication in $A$ and actions of $A$ on $N$, and $N^2 = 0$, the \textit{split null extension of $A$ by $N.$} The bimodule $N$ is called Jordan if $E$ is a Jordan algebra. Suppose that $N$ is Jordan, and let $e \in A$ be an idempotent. Then we may consider \textit{the Peirce decomposition of $N$ with respect to $e$}: $N = N_0 \oplus N_1 \oplus N_2$ defined by $N_i = N \cap E_i, i = 0, 1, 2.$ If $e$ is the unit of $A$, then $N$ is called \textit{unital} if $N = N_2(1)$ (in other words, $1$ acts identically on $N$).

The fact that $M^2 = 0$ allows us to apply representation theory to our problem:

\begin{Lem}
\label{ext_crit}
Let $J, \bar{J}, M, d, \bar{d}$ be as above. Then $M$ has a structure of a unital Jordan bimodule over $\bar{J}$ such that $\bar{d}(\bar{J})\cdot M = 0.$ On the other hand, let $\bar{J}$ be a simple Jordan algebra, $\bar{d}$ be a derivation with invertible values of $\bar{J}$ and $M$ be a unital $\bar{J}$-bimodule such that $\bar{d}(\bar{J})\cdot M = 0$. Consider the split null extension $E$ of $\bar{J}$ by $M$. Then the mapping $d'$ which extends $\bar{d}$ by $d'(M) = 0$ is a derivation with invertible values of $E.$
\end{Lem}
\textbf{Proof.} 
Denote the canonical projection alongside $M$ by $\beta.$ There exists a linear mapping $\delta: \bar{J} \rightarrow J$ such that $\beta\delta = 1_{\bar{J}}.$ Now, since $M^2 = 0$, from \cite[p. 92]{Jac} it follows that $M$ is a Jordan bimodule over $\bar{J}$ with the action given by
\begin{equation}
\label{action}
a\cdot m = m \cdot a = \delta(a)m, a \in \bar{J}, m \in M,
\end{equation}
where the right-hand sides are the products in $J.$ It is obvious that $\delta(\bar{1}) = 1 + m$ for some $m \in M.$ Therefore, $M$ is a unital bimodule over $\bar{J}.$ Recall that $d(M) = 0$. Hence, for $a \in J$ we have 
\begin{equation}
\label{zero_der_mult}
0 = d(am) = d(a)m.
\end{equation}
Let $\bar{a} = \beta(a) \in \bar{J}.$ Since $\beta\delta = 1_{\bar{J}}$, by the definition of $\bar{d}$ we have $\delta(\bar{d}(\bar{a})) = d(a) + m'$, where $m' \in M.$ Now let $m \in M$. Then $\bar{d}(\bar{a})m = \text{(by (\ref{action}))} = \delta(\bar{d}(\bar{a}))m = (d(a) + m')m = 0,$ by (\ref{zero_der_mult}) and the fact that $M^2 = 0.$ We proved the first part of the lemma. The second part of the lemma can be easily proved by a direct computation. \hfill$\Box$\vspace{0.3cm}

The lemma above gives us necessary and sufficient conditions for existence of nontrivial extensions described in theorem \ref{first_th}.

Let $J, \bar{J}$ and $\delta$ be as above. It is known \cite[p. 94]{Jac} that for $a, b \in \bar{J}$ we have $\delta(a)\delta(b) = \delta(ab) + h(a,b),$ where $h(a,b) \in M$, $h$ is a bilinear map and
\begin{equation}
\label{cohom_1}
h(a,b) = h(b,a),
\end{equation}
\begin{equation}
\label{cohom_2}
(h(a,a)\cdot b)\cdot a + h(a^2,b)\cdot a + h(a^2b,a) = a^2\cdot h(b,a) + h(a,a)\cdot (ba) + h(a^2,ba) = 0.
\end{equation}
The set of all such mappings $h$ modulo certain equivalence relation (rendering equal all $h$ such that the corresponding maps $\delta$ are obtained from each other by adding a mapping $\mu: M \to M$) is exactly the second cohomology group of $\bar{J}$ with coefficients in $M$. Particularly, the equivalence class of  $h$ is zero if and only if $\delta$ is a sum of homomorphism of algebras $\bar{J}$ and $J$ and a mapping $\mu: M \to M$ (hence the extension of $\bar{J}$ by $M$ is split). 

Now let $J$ be a simple finite-dimensional Jordan algebra and $d$ be a derivation with invertible values of $J$. Further in this section, as lemma \ref{ext_crit} suggests, we seek unital $J$-bimodules $M$ such that
\begin{equation}
\label{annul}
d(J)\cdot M = 0.
\end{equation}
Note that unital representations of simple finite-dimensional Jordan algebras were described in \cite{Jac}. %It was proved that every bimodule over simple finite-dimensional Jordan algebra is completely reducible, therefore, it suffices to consider only irreducible bimodules.

\subsection{The case $J$ is a division algebra of Albert type}
Let $J$ be a division algebra of Albert type. By \cite{Jac_rep} $M$ is isomorphic to the sum of the regular $J$-bimodules. Since $J$ is a unital algebra, in this case there are no $J$-modules with the required property (\ref{annul}), which together with lemma \ref{ext_crit} implies the 

\begin{proposition}
\label{ext_alb}
Let $J$ and $M$ be as in the beginning of the section. Suppose that $J/M$ is an algebra of Albert type. Then $M = 0$ and $J$ is a simple Jordan algebra.
\end{proposition}

\subsection{The case $J = J(V,f)$}Let $J = J(V,f)$ be an algebra of symmetric nondegenerate bilinear form and $d$ its derivation with invertible values. The bimodules over $J(V,f)$ with the property (\ref{annul}) are described in the  

\begin{Lem}
\label{form_mod}
Let $J = J(V,f)$ be an algebra of symmetric nondegenerate bilinear form over its center $Z$, let $d$ be a derivation with invertible values of $J$, and let $M$ be a unital bimodule over $J$ such that $d(J)\cdot M = 0.$ Then $V\cdot M = 0.$
\end{Lem}
{\bf Proof.} From (\ref{dt}) it follows that $d(V) \subseteq V.$ Suppose that $d(V) = 0.$ Then for any $z \in Z, v \in V$ we have $0 = d(zv) = d(z)v,$ hence $d(Z) = 0$ and $d = 0$. Therefore we may assume that there exists a nonzero $v \in Im(d) \cap V.$ Then, by lemma  \ref{form_der_char}, $v^2 = f(v,v) \neq 0.$  If $f(v,v) = \alpha^2$ is a square in $\mathbb{F},$ then for $u = \frac{1}{2\alpha}v$ we have $f(u,u) = \frac{1}{4}.$ If $f(v,v)$ is not a square in $\mathbb{F}$, we can take a quadratic extension $\mathbb{F}' = \mathbb{F}[\sqrt{f(v,v)}]$ of $\mathbb{F}$ and consider $J' = J \otimes_F F',$ which is again an algebra of symmetric nondegenerate bilinear form. Also $J'$ acts on $M' = M \otimes_F F'$ in a natural way, and clearly $v \cdot M' = 0$. Therefore, we may assume that $f(v,v) = \frac{1}{4}.$ Note that the element $e = \frac{1}{2} + v$ is an idempotent of $J.$ The Peirce decomposition of $J$ with respect to $e$ is the following: $J_0 = \langle \frac{1}{2} - v \rangle, J_1 = \{ u \in V: f(u,v) = 0 \}, J_2 = \langle \frac{1}{2} + v \rangle.$ Since $v\cdot M = 0$, for any $m \in M$ we have $e\cdot m = \frac{1}{2}m,$ hence $M = M_1(e).$ Thus $J_1\cdot M \subseteq \text{(by (\ref{pd}))} \subseteq M_0(e) + M_2(e)= 0,$ and $V \cdot M = (\langle v\rangle + J_1)\cdot M = 0.$ \hfill$\Box$\vspace{0.3cm}

\medskip

The lemma implies that the action on this bimodule is the action given by 
\begin{equation}
\label{form_action}
1\cdot m = m, v \cdot m = 0 \text{ for }v \in V, m \in M.
\end{equation} This bimodule is indeed Jordan: note that the split null extension of $J$ by $M$ is equal to $J(V \oplus M,f')$, where $f'$ extends $f$ such that $f'(M,V \oplus M) = 0.$ We now show that extensions by bimodules of this type always split:

\begin{Lem}
Let $A$ be an extension of an algebra $J$ of nondegenerate bilinear form $J(V,f)$ by an ideal $M$ such that $M^2 = 0.$ Suppose that the actions of $J$ on $M$ given by (\ref{action}) and (\ref{form_action}) coincide. Then $A$ is the split null extension of $J$ by $M.$
\end{Lem}
\textbf{Proof.} 
 To prove that $A$ is split we have to show that $J$ is a subalgebra of $A$. Hence the map $\delta$ must be a sum of homomorphism of algebras and a mapping $\mu: M \to M$, or equally $h$ must be 0. It is easy to see that subtracting a mapping $\mu: M \to M$ we can choose $\delta$ such that $\delta(1) = 1,$ the unit of $A.$ Hence for any $a \in J, \delta(a)\delta(1) = \delta(a),$ therefore $h(F,J) = h(J,F) = 0.$ Now, let $a, b \in V$ in (\ref{cohom_2}). Using the fact that $ab, a^2 \in F$ and $V\cdot M = 0$, we have $h(a,a)\cdot(ba) =  h(a,a)f(a,b) = 0.$ Since $f$ is nondegenerate, we can choose $b$ such that $f(a,b) \neq 0,$ thus $h(a,a) = 0.$ Therefore, (\ref{cohom_1}) implies that $h(V,V) = 0$ and $h= 0,$ hence $\delta$ is an algebra homomorphism and $A$ is split. \hfill$\Box$\vspace{0.3cm}

\medskip

Our considerations can be summed up in the

\begin{proposition}
\label{ext_form}
Let $J$ and $M$ be as in the beginning of the section. Suppose that $J/M$ is an algebra of nondegenerate bilinear form. Then $J$ is an algebra of bilinear form, and $M$ is the radical of $J.$
\end{proposition}

It is not hard to construct derivations with invertible values for algebras $J(V,f)$ where $f$ may be degenerate, but not zero (use, for example, lemma \ref{form_der}).

\subsection{The case $J = A^{(+)}$ or $H(A,*), A$ is simple associative}
In this section we consider the cases where $J$ is either $A^{(+)}$, or $H(A,*)$, where $A = D$ or $D_2$.   

Note that in each case $J$ is a central simple finite-dimensional special algebra over $Z = Z(D)$. Let $\bar{Z}$ be the algebraic closure of $Z$, and consider the scalar extension $\bar{J}$ of $J$ by $\bar{Z}$. Clearly $\bar{J}$ is a central simple special Jordan algebra over an algebraically closed field. If degree of $J$ is $\leq 2$ then the classification of central simple finite-dimensional Jordan algebras over an algebraically closed field implies that  $\bar{J}$ is either $\bar{Z}$ itself or an algebra of nondegenerate symmetric bilinear form \cite[p. 204]{Jac}. Note that $J$ is special, Hence, if degree of $J$ is $\geq 3$, from the same classification it follows that  $\bar{J}=M_n(\bar{Z})^{(+)}$, or $H(M_n(\bar{Z}),tr)$ or $H(M_{2n}(\bar{Z}),Symp)$, where $n \geq 3$.

We will prove that $J$ has no nonzero bimodules with the property (\ref{annul}). Moreover, we will prove that any $J$-bimodule $M$ is faithful, that is, there is no element $j \in J$ such that $j \cdot M = 0$. Clearly it suffices to prove the faithfulness of all $\bar{J}$-bimodules. Hence, we have to study bimodules over these algebras. Note that we have already studied the case where $\bar{J}$ is an algebra of a symmetric bilinear nondegenerate form (the case $\bar{J} = \bar{Z}$ is obvious). Therefore we may assume that the degree of $J$ is $\geq 3.$ Thus, we have to prove the faithfulness of bimodules over algebras $M_n(\bar{Z})^{(+)}, H(M_n(\bar{Z}),T)$ and $H(M_{2n}(\bar{Z}),Symp)$, where $n \geq 3$. Since bimodules over these algebras are completely reducible \cite[p. 273]{Jac}, it suffices to study only irreducible ones.

Let us list here the irreducible bimodules over these algebras :\\
1) Let $J = H(A,*)$, where $(A,*) = (M_n,T)$ or $(M_{2n},Symp)$. Then every irreducible unital bimodule over $J$ is isomorphic to $H(A,*)$ or $S(A,*)$ (the set of $*$-skew elements) with the action $h\cdot x = \frac{1}{2}(xh + hx)$.\\
2) Let $J = M_n^{(+)}$. Then every irreducible unital bimodule over $J$ is isomorphic to one of the following: (i) $M_n$, with the action of $b \in M_n^{(+)}$ given by $b \cdot x = \frac{1}{2}(xb + bx),$ (ii) $H(M_n,T)$ (the space of symmetric $n \times n$-matrices), with the action of $b \in M_n^{(+)}$ given by $b \cdot x = \frac{1}{2}(xb^T + bx),$ (iii) $H(M_n,T)$ with $b \cdot x = \frac{1}{2}(xb + b^Tx),$ (iv) $S(M_n,T)$ (the space of skew-symmetric $n \times n$-matrices) with $b \cdot x = \frac{1}{2}(xb^T + bx),$ (v) $S(M_n,T)$ with $b \cdot x = \frac{1}{2}(xb + b^Tx).$

\begin{Lem}
Let $n \geq 3.$ Then bimodules listed of the type $1)$ are faithful.
\end{Lem}
\textbf{Proof.} It is clear that in both cases $H(A,*)$ is a faithful module over itself (because it is a unital algebra). The case $S(A,*)$ can also be done somewhat independently of the choise of $(A,*).$ Note that the condition $x\cdot S(A,*) = 0$ for some $x \in H(A,*)$ means that $x$ anticommutes with all $S(A,*)$ in the enveloping associative algebra. It is an easy exercise in linear algebra that in both cases every element of $A$ can be written in the form $s_{11}s_{12} + \ldots + s_{k1}s_{k2},$ where $s_{ij} \in S(A,*)$ for all $i, j$. Hence $x$ lies in the center of associative algebra $A$, which is the ground field in both cases. But it is obvious that the elements of the ground field cannot act as zero operators on $S(A,*).$ \hfill$\Box$\vspace{0.3cm}

\begin{Lem}
Let $n \geq 3.$ Then bimodules of the type $2)$ are faithful.
\end{Lem}
\textbf{Proof.} The bimoudle (i) is clearly faithful, since $M_n^{(+)}$ is a unital algebra. Note that it suffices to prove that modules (ii) and (iv) are faithful, since actions on modules (iii) and (v) can be obtained using corresponding actions on (ii) and (iv) and transpositions. Proving the faithfulness of modules (ii) and (iv) is entirely a matter of matrix algebra: it suffices to write the action of an arbitrary matrix $b \in M_n$ on the bases of $Sym_n$ and $Skew_n$ correspondingly and verify that no nonzero matrix can nullify all basis elements. \hfill$\Box$\vspace{0.3cm}

We can sum up our considerations in the following

\begin{proposition}
\label{ext_deggeq3}
Let $J$ and $M$ be as in the beginning of the section. Suppose that $J/M$ is an algebra $H(A,*)$ and that it has degree $\geq 3$. Then $M = 0$ and $J$ is a simple Jordan algebra.
\end{proposition}

Finally, the results of propositions \ref{ext_alb}, \ref{ext_form} and \ref{ext_deggeq3} imply the main result of the section:

\begin{Th}
Let $J$ be a finite-dimensional Jordan algebra admitting a derivation with invertible values. Then $J$ is either simple or an algebra of a symmetric bilinear form (possibly degenerate).
\end{Th}

Note that the propositions \ref{ext_alb} and \ref{ext_form} were obtained without imposing finite-dimensionality of any space involved. We only used finite-dimensionality over center while studying the algebras $A^{(+)}$ and $H(A,*)$, in fact implicitly using the classification of simple finite-dimensional associaive algebras over an algebraically closed field. In general, infinite-dimensional case we would have to study the representations of algebras $A^{(+)}$ and $H(A,*)$ for $A$ a (possibly infinite-dimensional over its center) division algebra $D$ or $2 \times 2$ algebra over $D,$ which does not appear feasible. Therefore, this approach cannot be extended any further.

\section{Appendix}
In the paper \cite{KP} the first and the third author described derivations with invertible values of algebras $C$, which are split Cayley-Dickson algebras over its center $Z$. During the consideration it was assumed that a derivation $d$ of $C$ with invertible values is $Z$-linear, but this was given without proof. Here we give the proof of this fact:

\begin{Lem}
Let $C$ be a split Cayley-Dickson algebra over its center $Z,$ and let $d$ be a derivation with invertible values of $C.$ Assume that the characteristic of the base field $\mathbb{F}$ is not 2. Then $d$ is $Z$-linear.
\end{Lem}
\textbf{Proof.} It is well known that the algebra $C^{(+)}$ is an algebra of symmetric nondegenerate bilinear form of dimension 8 over $Z$. Moreover, it is easy to find an isotropic subspace $V$ of $C$ of $Z$-dimension 4. Indeed, using the matrix representation of elements of $C$, one can take 

$$V = \{\begin{pmatrix} \alpha & v \\ 0 & 0 \end{pmatrix}, \alpha \in Z, v \in Z^3 \}.$$

Then by lemma \ref{form_der_lin}, $d$ is linear over $Z$. \hfill$\Box$\vspace{0.3cm}

%{\bf Acknowledgements.}
%The authors are grateful to Prof. Victor Zhelyabin (Sobolev Institute of Math., Russia) for noting the gap in the article \cite{KP}.

\medskip

%\newpage

\end{document}